\documentclass[12pt]{article}
\usepackage{amsmath}
\usepackage{amssymb}
\usepackage{amsthm}
\usepackage{url}
\usepackage{enumerate}
\usepackage{amsfonts,amsmath,amssymb,graphicx,psfrag}
\usepackage[T1]{fontenc}
\usepackage[utf8]{inputenc}
\usepackage{mathrsfs,amsmath}
\numberwithin{equation}{section}
\newtheorem{theorem}{Theorem}[section]
\newtheorem{proposition}[theorem]{Proposition}

\newtheorem{Remark}[theorem]{Remark}

\newcommand\dstyle\displaystyle

\newcommand\al\alpha
\newcommand\be\beta
\newcommand\de\delta
\newcommand\la\lambda
\newcommand\tha\theta

\newcommand\iy\infty

\newcommand{\hyp}[5]{\,\mbox{}_{#1}F_{#2}\!\left(
  \genfrac{}{}{0pt}{}{#3}{#4};#5\right)}

\newcommand\bma{\begin{pmatrix}}
\newcommand\ema{\end{pmatrix}}

\begin{document}

\title{Meixner $d$-Orthogonal Polynomials Arising from $\mathfrak{su}(1,1)$ }
\author{B. Halouani and   F.  Bouzeffour. \\ Department of mathematics, College of Sciences.\\
 King Saud University, P. O Box $2455$ Riyadh $11451$, Saud Arabia. \\
halouani@ksu.edu.sa, fbouzaffour@ksu.edu.sa}

\maketitle
\begin{abstract}In this study, we present a novel family of Meixner-type $d$-orthogonal polynomials, which are distinguished as a particular subset of multiple orthogonal polynomials. We demonstrate their connection to the Lie algebra $\mathfrak{su}(1,1)$ by identifying them as matrix elements of an appropriately defined nonlinear operator. Utilizing Barut-Girardello coherent states, we explicitly outline their key features, including recurrence relations, generating functions, and $d$-orthogonality relations, among others.
\end{abstract}

\section{Introduction}
The synergy between Lie theory and orthogonal polynomials is well-recognized, highlighted by pivotal references such as \cite{M2,G3,Bou,Bouz,Bouze}. These polynomials notably emerge as matrix elements associated with the  generators, paving the way for an algebraic framework that facilitates the derivation of generating functions, recurrence relations, and orthogonality properties of these polynomials.

This paper aims to further elucidate this framework by exploring the relationship between the Lie algebra $\mathfrak{su}(1,1)$ and $d$-orthogonal polynomials. The latter are a subset of multiple orthogonal polynomials, recognized for their connections to solutions of some recurrence relations, as demonstrated in \cite{A2,A3,Bou,Bouz,Bouze} and further supported by \cite{A1,B1,I1,G1,G2,G3,G4,K1,O1,P1,S1}.

Our investigation focuses on the three-dimensional Lie algebra $\mathfrak{su}(1,1)$, defined by the generators $J_0$, $J_{-}$, and $J_{+}$. We analyze the matrix elements stemming from the operator
\begin{equation*}
S=e^{ J_+}e^{Q(J_-)}.
\end{equation*}
The study of this operator not only extends the classical Meixner polynomials but also uncovers essential features such as generating function, recurrence relations, and differential equations through the lens of Barut-Girardello coherent states, as referenced in \cite{B1,I1,S1}.

The paper is structured as follows: Section 2 revisits the foundational concepts and key findings regarding $d$-orthogonal polynomials. Section 3 introduces the $\mathfrak{su}(1,1)$ algebra and its representations, detailing the Barut-Girardello coherent states and establishing significant identities within $\mathfrak{su}(1,1)$. Section 4 presents an in-depth analysis of the operator $S$ and its related matrix elements expressed through $d$-OPS. This algebraic approach allows us to derive crucial properties, such as biorthogonality relations, recurrence formulas, generating functions, and the application of the lowering operator. Finally, Section 5 discusses our main results concerning a new family of Meixner-type $d$-OPS, explicitly defining a linear functional vector that confirms their $d$-orthogonality via biorthogonality relations.

\section{$d$-Orthogonal Polynomials}
Consider $\mathcal Q$ as the linear space composed of polynomials with complex coefficients, and $\mathcal Q'$ as its algebraic dual. The action of a linear functional $u \in \mathcal Q'$ on a polynomial $f \in \mathcal Q$ is denoted by $\langle u,f \rangle$. A sequence of polynomials $\{Q_n\}_{n\geq0}$ is termed a polynomial set if deg$(Q_n) = n$ for every nonnegative integer $n$.

The concept of $d$-orthogonality extends the traditional notion of usual orthogonality, as introduced by Van Iseghem and Maroni (refer to \cite{I1,M1}). Let $\{Q_n\}_n$ be a polynomial set within $\mathcal Q$ and consider $d$ as a positive integer. The set $\{Q_n\}_{n\geq0}$ is defined as a $d$-orthogonal polynomial set ($d$-OPS) relative to the $d$-dimensional functional vector $\overrightarrow{\mathcal L} = (\mathcal L_0, \mathcal L_1, \ldots, \mathcal L_{d-1})$ if, for each integer $i \in \{0,1,\ldots,d-1\}$,
\begin{align*}
&\langle \mathcal L_i,Q_mQ_n\rangle = 0, \quad \text{for } n \geq md + i + 1,\\
&\langle \mathcal L_i,Q_nQ_{nd+i}\rangle \neq 0, \quad \text{for } n \geq 0.
\end{align*}
This definition reverts to the familiar concept of orthogonality when $d=1$.

The study and development of $d$-orthogonal polynomials, along with multiple orthogonal polynomials, have captivated many researchers, offering numerous explicit examples and extensive theoretical frameworks (see, for example, \cite{B5}).

Furthermore, a polynomial sequence $\{Q_n\}_{n\geq0}$ is recognized as $d$-OPS if it adheres to a recurrence relation of order $d+1$, described by
\begin{equation}
xQ_n(x) = \beta_{n+1}Q_{n+1}(x) - \sum_{i=0}^d \alpha_{n,i} Q_{n-i}(x), \label{eq0}
\end{equation}
where $\beta_{n+1}\alpha_{n,d} \neq 0$, under the convention that $Q_{-n} = 0$ for $n\geq1$. For $d=1$, this scenario simplifies to what is known as the Favard Theorem.
\subsection{Non linear coherent states.}
A non linear coherent states is an {\em overcomplete} set of vectors
in a  Hilbert space, labelled by a continuous parameter $z$ which
runs over a complex domain. The vectors are, in addition, subject to
a {\em resolution of the identity condition\/.} More precisely, let
$H$ be a (complex, separable, infinite dimensional) Hilbert space,
with inner product $\langle\ |\ \rangle$, $\{e_n\}_{n=0}^\infty$ an
orthonormal basis of it and let $\{x_n\}_{n=0}^\infty, \; x_0 = 0$,
be an infinite sequence of positive numbers. Let $\lim_{n\rightarrow
\infty}x_n = L^2$, where $L > 0$ could be finite or infinite, but
not zero. We shall use the notation $x_n! = x_1 x_2 \cdots x_n$ and
$x_0! = 1$. For each $z\in \mathcal D$ (some domain in $\mathbb C$),
we define a non--linear coherent state, i.e., a vector $v_z \in H$,
in the manner
 \begin{equation}
  v_z = \mathcal N(\vert z\vert^2 )^{-\frac 12} \sum_{n=0}^\infty\frac {z^n}{\sqrt{x_n!}}e_n\;
\label{CS-def}
\end{equation}
where the normalization constant $\mathcal N(\vert z\vert^2 ) =
\sum_{n=0}^\infty\dfrac {\vert z\vert^{2n}}{x_n!}$ is chosen so that
$\Vert v_z \Vert = 1$. It is clear that the vectors $v_z$ are well
defined for all $z$ for which the above sum, representing $\mathcal
N(\vert z\vert^2 )$, converges, i.e., $\mathcal D = \{ z\in \mathbb
C \mid \vert z \vert < L\}$. Furthermore, we require that there
exist a measure $d\nu (z, \overline{z})$ on $\mathcal D$ for which
the resolution of the identity condition,
\begin{equation}
 \int_{\mathcal D} \vert v_z\rangle\langle v_z \vert\;\mathcal N(\vert z\vert^2 )\; d\nu (z, \overline{z}) =
     I_H\; ,
\label{resolid}
\end{equation}
holds. The symbol $P_z=\vert v_z\rangle\langle v_z \vert$ stands for
the projection operator onto the coherent state $v_z$. That is, for
any vector $u\in H$ we have $P_zu=\langle v_z|u\rangle v_z.$

It is easily seen that in order for (\ref{resolid}) to be satisfied,
$d\nu$ has to have the form, \begin{equation}
  d\nu (z , \overline{z}) = \frac 1{2\pi}\; d\theta\; d\lambda(r)\; , \qquad z = re^{i\theta}\;
\label{orthog-meas}
\end{equation}
where the measure $d\lambda$ is a solution of the moment problem,
\begin{equation}
\int_0^L r^{2n} d\lambda(r) = x_n! , \qquad n =0,1,2, \ldots  \;
\label{mom-prob}
\end{equation}
provided such a solution exists.  In most of the cases that occur in
practice, the support of the measure $d\nu$ is the whole of
$\mathcal D$, i.e., $d\lambda$ is supported on the entire interval
$(0, L)$.\\
Below are some examples of the above general construction.
\subsection{$\mathfrak{su}(1,1)$ Lie algebra.}

The Lie algebra $\mathfrak{su}(1,1)$ is the $3$-dimensional vector space with the basis elements $J_0,\, J_{-},\,J_{+}$ and
commutation relations
\begin{equation}\label{eq1}
[J_-,J_+]=2J_0,\,\,\,\,[J_0,J_{\pm}]=\pm J_{\pm}.
\end{equation}
The invariant element $C=J_+J_-+J^2_0-J_0$ for any representation
is $C=-\lambda(\lambda-1)I$, where $I$ is the identity operator. Thus a representation of $\mathfrak{su}(1,1)$ is determined by $\lambda$.\\
We are mainly interested by the positive discrete series $D_\beta$, where
 $\beta>0$. The corresponding Hilbert space is spanned by the orthonormal basis $|n\rangle\,\,n=0,1,\dots$.\\
  The action of the generators $J_0$ and $J_{\pm}$ is defined
as follows
\begin{eqnarray} J_+|n\rangle&=&\sqrt{(n+1)(n+\beta)}|n+1\rangle,
\\ J_-|n\rangle&=&\sqrt{n(\beta+n-1)}|n-1\rangle,\ \ n=1,2,\ldots, \ J_-e_0=0,\\
 J_0|n\rangle&=&(n+\frac{\beta}{2})|n\rangle.
\end{eqnarray}\label{K0}
For this representation, we have the conjugation relations
\begin{equation} J^\ast_0=J_0\,\,\,\text{
and}\,\,\, J^\ast_{\pm}=J_{\mp}.\label{St} \end{equation}  All the vectors $|n\rangle$,
$n=0,1,\dots$ can be obtained from $|0\rangle$ by acting by powers of the
raising operator $J_+$ and the vector $|n-i\rangle$ is also obtained from
$|n\rangle$ $(0\leq\,i\,\leq\,n)$ by acting by  powers of the lowering operator $J_-$; e.g.
\begin{align}
&J^n_+|0\rangle=\sqrt{n!(\beta)_n}|n\rangle,\\
&J^i_-|n\rangle=\sqrt{\frac{n!(\beta)_n}{(n-i)!(\beta)_{n-i}}}|n-i\rangle,
\end{align}
where we have used the common notation for the Pochhammer symbol
\begin{equation}
(a)_n=a(a+1)\dots (a+n-1),\ \ \ (a)_0=1.
\end{equation}

Let $x_n = \sqrt{n(n+\beta-1)}$, so that $L = \infty$. In that case the
coherent states $v_z$ is defined for all $z\in \mathbb C$
\begin{equation}
  |z\rangle = \sum_{n=0}^\infty\frac {z^n}{\sqrt{n! (\beta)_n}}\; |n\rangle \; , \qquad z \in \mathbb C.
\label{bar-gir-cs}
\end{equation}
According to (2.10) and (2.13) we get
\begin{equation}
 |z\rangle=\hyp01{-} {\beta}{zJ_+}|0\rangle
\end{equation}
where ${ }_rF_s$ is the hypergeometric function defined by
$$\hyp rs{a_1,a_2,\ldots,a_r}{b_1,b_2,\ldots,b_s}{x}:=\sum_{m=0}^\infty\frac{(a)_m(a_2)_m\cdots(a_r)_m}{(b_1)_m(b_2)_m\cdots(b_s)_m}\frac{x^m}{m!}$$

The normalized vector $|v_z\rangle=\frac{|z\rangle}{\||z\rangle\|}$ coincides with the
well--known the Barut--Girardello coherent states (BGCS), (see[4]).\\ These coherent states satisfy the resolution of the
identity,
\begin{equation}
   \frac 2\pi\; \int_{\mathbb C}\mid v_z\rangle
   \langle v_z \mid\; K_{\beta -1}(2r)\; I_{\beta -1}
          (2r )\; r\; dr\; d\theta = I\; , \qquad z = r e^{i\theta}\; ,
\end{equation}
where again, $K_\nu (x)$ is the order-$\nu$ modified Bessel function
of the second kind.\\The BGCS are eigenvalue of the lowering
operator $J_-$, i.e
\begin{equation}
J_-|z\rangle=z|z\rangle.
\end{equation}
For an analytic function
\begin{equation}\label{eq2}
F(J_-)|z\rangle=F(z)|z\rangle
\end{equation}
The Peremolove's coherent states (PCS) are defined for all $z$ in the
unit disc $\mathcal D = \vert z\vert < 1$, by
\begin{align}
|z\rangle = S(\xi)|0\rangle =(1-|z|^2)^{\beta/2} \sum_{n=0}^\infty \sqrt{\frac{(\beta)_n}{n!}}z^n|n\rangle,
\end{align}
where $\xi=\tanh^{-1}|z|e^{i\theta}$ and $S(\xi)=e^{\xi
J_+-\overline{\xi}J_-}$ is the $su(1,1)$ displacement operator
\cite{P1}.\\The corresponding moment problem is now \begin{equation}
  \int_0^1 r^{2n} \; d\lambda (r) =  \frac {n!}{(\beta)_n}\; ,
\label{su-mom-prob}
\end{equation}
which has the solution, \begin{equation}
  d\lambda (r) = 2(\beta -1)\; r(1 - r^2)^{\beta -2}\; dr, \qquad 0 \leq r < 1\; .
\label{su-meas}
\end{equation}
\subsection{Some identities in $\mathfrak{su}(1,1)$}
We denote by $\mathcal{U}(\mathfrak{su}(1,1))$ the universal enveloping algebra of $\mathfrak{su}(1,1)$. Utilizing the relations (\ref{eq1}), we can establish through induction on $n$ that
\begin{align}
[J_0,J^n_{\pm}] &= \pm nJ^{n}_{\pm}, \label{identity1} \\
[J_-^n,J_+] &= 2nJ_0J^{n-1}_- + n(n-1)J^{n-1}_-. \label{identity2}
\end{align}
From these identities  \ref{identity1}) and \eqref{identity2} , and for any polynomial $Q$, it follows that
\begin{align}
[J_0 ,Q(J_\pm)] &= \pm J_\pm Q'(J_\pm), \label{identity3} \\
[Q(J_-),J_+] &= 2J_0Q^{'}(J_-) + J_-Q^{''}(J_-). \label{identity4}
\end{align}
These results, in conjunction with the Baker--Hausdorff formula (see for instance, \cite[\S 2.15]{V1}),
\begin{equation}
e^ABe^{-A} = B + [A,B] + \frac{1}{2!}[A,[A,B]] + \frac{1}{3!}[A,[A,[A,B]]] + \dots, \label{eqBakerHausdorff}
\end{equation}
yield the following relations:
\begin{align}
e^{Q(J_{\pm})}J_0e^{-Q(J_{\pm})} &= J_0 \mp J_{\pm}Q^{'}(J_{\pm}), \label{relation1} \\
e^{-Q(J_-)}J_+e^{Q(J_-)} &= J_+ - 2J_0Q^{'}(J_-) + J_-(Q^{'}(J_-)^2 - Q^{''}(J_-)). \label{relation2}
\end{align}

\section{Matrix Elements of an Operator and $d$-Orthogonal Polynomials}
Let $r \geq 1$ and $d = 2r +1$ be positive integers, and let $Q$ be a polynomial of degree $r$ with $Q(0) = 0$. We define the operator $S$ as
\begin{equation}
S = e^{J_+}e^{Q(J_-)}. \label{eq12}
\end{equation}
The operator $S$ and its matrix elements will be the focal points of our study for the remainder of this paper. It is important to note that $S$ is invertible, with its inverse given by $S^{-1} = e^{-Q(J_-)}e^{-J_+}$. The matrix elements of $S$ and $S^{-1}$ are denoted and defined by
\begin{equation}\label{br1} \psi_{nk} = \langle \beta,k \lvert S \rvert n,\beta \rangle, \quad \phi_{nk} = \langle \beta,n \lvert S^{-1} \rvert k,\beta \rangle. \end{equation}
These matrix elements $\psi_{nk}$ and $\phi_{nk}$ satisfy the biorthogonality relation
\begin{equation}
\sum_{k=0}^{\infty} \psi_{nk} \phi_{mk} = \delta_{nm}. \label{eq13}
\end{equation}
To demonstrate this, consider that
\[ S \lvert n,\beta \rangle = \sum_{k=0}^{\infty} \langle \beta,k \lvert S \rvert n,\beta \rangle\,\, \lvert k,\beta \rangle, \]
leading to
\[ S^{-1}S \lvert n,\beta \rangle = \sum_{k=0}^{\infty} \sum_{r=0}^{\infty} \langle \beta,k \lvert S \rvert n,\beta \rangle \langle \beta,r \lvert S^{-1} \rvert k,\beta \rangle \lvert r,\beta \rangle. \]
Then \begin{align*} \langle \beta,m \lvert S^{-1}S \rvert n,\beta \rangle &= \sum_{k=0}^{\infty}\sum_{r=0}^{\infty} \langle \beta,k \lvert S \rvert n,\beta \rangle \langle \beta,r \lvert S^{-1} \rvert k,\beta\rangle\langle \beta,m \lvert r,\beta \rangle\\&=\sum_{k=0}^{\infty} \langle \beta,k \lvert S \rvert n,\beta \rangle \langle \beta,m \lvert S^{-1}  \lvert k,\beta \rangle. \end{align*}
Therefore
\[ \sum_{k=0}^{\infty} \langle \beta,k \lvert S \rvert n,\beta \rangle \langle \beta,m \lvert S^{-1} \rvert k,\beta \rangle=\delta_{n,m}. \]
The results follows from \eqref{br1}.
\subsection{Recurrence Relation}
To establish a recurrence relation satisfied by $\psi_{nk}$, we consider the expression $\langle \beta,n \lvert J_0S \rvert n,\beta \rangle$. By \ref{eq1}, we have
\begin{equation}
\langle \beta,k \lvert J_0S \rvert n,\beta \rangle = \left(k+\frac{\beta}{2}\right)\langle \beta,k \lvert S \rvert n,\beta \rangle = \left(k+\frac{\beta}{2}\right)\psi_{nk}. \label{eq14}
\end{equation}
On the other hand,
\begin{equation}
\langle \beta,k \lvert J_0S \rvert n,\beta \rangle = \langle \beta,k \lvert S(S^{-1}J_0S) \rvert n,\beta \rangle. \label{eq15}
\end{equation}
According to \eqref{relation1} and \eqref{relation2}, we have
\begin{align*}
S^{-1}J_0S &= e^{-Q(J_-)}\left(e^{-J_+}J_0e^{J_+}\right)e^{Q(J_-)} \\
&= e^{-Q(J_-)}\left(J_+ + J_0\right)e^{Q(J_-)} \\
&= e^{-Q(J_-)}J_+e^{Q(J_-)} + e^{-Q(J_-)}J_0e^{Q(J_-)} \\
&= J_+ + J_0\left(1-2Q'(J_-)\right) + J_-\left(Q'(J_-)^2 + Q'(J_-) - Q''(J_-)\right).
\end{align*}
Since $\deg(Q) = r$, then we can write
\begin{align*}
1-2Q'(x) &= \sum_{i=0}^{d} a_ix^i, \\
x\left(Q'(x)^2 + Q'(x) - Q''(x)\right) &= \sum_{i=0}^d b_ix^i,
\end{align*}
where $d=2r+1,$ $a_i$ and $b_i$ are complex numbers with $a_i = 0$ for $r \leq i \leq d,\ a_{r-1}b_d \neq 0$.
Therefore, (\ref{eq15}) becomes
\begin{align*}
\langle \beta,k \lvert J_0S \rvert n,\beta \rangle &= \langle \beta,k \lvert SJ_+ \rvert n,\beta \rangle + \langle \beta,k \lvert SJ_0\sum_{i=0}^d a_iJ_-^i \rvert n,\beta \rangle + \langle \beta,k \lvert S\sum_{i=0}^d b_iJ_-^i \rvert n,\beta \rangle.
\end{align*}
From \ref{eq1}, \eqref{identity1} and  \eqref{identity2}, we obtain
\begin{equation}
\langle \beta,k \lvert J_0S \rvert n,\beta \rangle = \sqrt{(n+1)(n+\beta)}\psi_{n+1,k} + \sum_{i=0}^{d} \gamma_{ni}\psi_{n-i,k},
\end{equation}
where
\[\gamma_{ni} = \left(a_i\left(n-i+\frac{\beta}{2}\right)+b_i\right)\sqrt{\frac{n!(\beta)_n}{(n-i)!(\beta)_{n-i}}}.\]
Comparing this with (\ref{eq14}) leads to the recurrence relation:
\begin{equation*}
\left(k+\frac{\beta}{2}\right)\psi_{nk} = \sqrt{(n+1)(n+\beta)}\psi_{n+1,k} + \sum_{i=0}^{d} \gamma_{ni}\psi_{n-i,k}.
\end{equation*}
By dividing each side of the above equation by $\psi_{0k}$, we derive the following proposition:

\begin{proposition} The matrix elements $\psi_{nk}$ are expressed as
\[\psi_{nk} = P_n(k)\psi_{0k},\]
where $\{P_n(k)\}_{n\geq0}$ is a $d$-OPS in the argument $k$, satisfying the
recurrence relation of order $d+1$
\begin{equation}
\left(k + \frac{\beta}{2}\right)P_n(k) = \sqrt{(n+1)(n+\beta)}P_{n+1}(k) + \sum_{i=0}^{d} \gamma_{ni}P_{n-i}(k), \label{eq19}
\end{equation}
with the initial conditions $P_0(k) = 1$ and $P_n(k) = 0$ for $n < 0$.
\end{proposition}

\subsection{Generating function}
To facilitate the analysis of our operator \(S\) and its impact on the \(d\)-orthogonal polynomial sequence, we will need to work with the associated monic polynomials \(\widehat{P}_n(k)\). These polynomials are explicitly defined by the formula:
\begin{equation}
\widehat{P}_n(k) = \sqrt{n!(\beta)_n}P_n(k), \label{eq20Refined}
\end{equation}
where $P_n(k)$ denotes the $n$-th polynomial in the sequence \eqref{eq19}.

Turning our attention to the operator \(S\), we examine its effect on the vectors \(\lvert z, \beta \rangle\). This analysis leads us to the generating function \(F(z, k)\), which encapsulates the properties of the normalized polynomials \(\{\widehat{P}_n(k)\}_{n\geq0}\). This function provides a comprehensive view of the polynomial sequence through its representation as an infinite sum:
\begin{equation}
F(z, k) := \sum_{n=0}^\infty \frac{\widehat{P}_n(k)}{n!(\beta)_n} z^n,
\end{equation}
where each term of the series integrates the monic polynomial \(\widehat{P}_n(k)\), normalized by the factorial \(n!(\beta)_n\), and scaled by the power \(z^n\).

We have
\begin{equation}
\langle \beta,k|S|z,\beta\rangle= \langle \beta,k|S| \sum_{n=0}^\infty\frac {z^n}{\sqrt{n! (\beta)_n}} |n,\beta\rangle
=\sum_{n=0}^\infty\frac{\psi_{nk}}{\sqrt{n! (\beta)_n}}z^n\label{eq21}
\end{equation}
and
\begin{equation}
\langle \beta,n|S|z,\beta\rangle= \langle \beta,n|e^{J_+}e^{Q(J_-)}|z,\beta\rangle.\label{eq22}
\end{equation}
Taking into account of equation \eqref{eq2} and  the following identity
$$e^t\hyp01{-} {\beta}{zt} =\sum_{m=0}^\infty\hyp11{-m}{\beta}{-z}\frac{t^m}{m!},$$
we obtain
\begin{align*}
\langle \beta,n|S|z,\beta\rangle&=e^{Q(z)}\langle \beta,n|e^{J_+}\hyp01{-} {\beta}{zJ_+}|0\rangle\\
&=e^{Q(z)}\sum_{m=0}^\infty\frac{1}{m!}\hyp11{-m}{\beta}{-z}\langle \beta,n|J_+^m|0\rangle\\
&=e^{Q(z)}\sum_{m=0}^\infty\sqrt{\frac{(\beta)_m}{m!}}\hyp11{-m}{\beta}{-z}\langle \beta,n|m,\beta\rangle\\
&=\sqrt{\frac{(\beta)_k}{k!}}e^{Q(z)}\hyp11{-k}{\beta}{-z}.
\end{align*}
It follows from (\ref{eq21}) that the matrix elements $\psi_{nk}$ are generated by
\begin{equation}
\sum_{n=0}^\infty\frac{\psi_{nk}}{\sqrt{n! (\beta)_n}}z^n=\sqrt{\frac{(\beta)_k}{k!}}e^{Q(z)}\hyp11{-k}{\beta}{-z}.\label{eq23}
\end{equation}
Since $Q(0)=0$, it is clear from (\ref{eq23}) that
\begin{equation*}
 \psi_{0k}=\sqrt{\frac{(\beta)_k}{k!}}.
\end{equation*}
Note that the left hand side of \eqref{eq23} is an analytic function in the variable $z$. Therefore, the series in \eqref{eq23} converges throughout the entire complex space $\mathbb{C}$. This demonstrates the following proposition:
\begin{proposition}The $d$-OPS $\{\widehat{P}_n(k)\}_{n\geq0}$ is generated by
\begin{equation}
\sum_{n=0}^\infty\frac{\widehat{P}_n(x)}{n! (\beta)_n}z^n=e^{Q(z)}\hyp11{-x}{\beta}{-z}.\label{eq25}
\end{equation}
\end{proposition}
\subsection{Connection with the First Kind Meixner Orthogonal Polynomials}
The Meixner polynomials of the first kind, denoted as $\mathcal{M}_n(x,\beta,c)$, are generated by the following series expansion (see \cite{K3,K4}):
\begin{equation}
\sum_{n=0}^\infty \mathcal{M}_n(x,\beta,c) \frac{z^n}{n!} = e^z \hyp11{-x}{\beta}{\frac{1-c}{c}z},\label{eq26}
\end{equation}
where $\beta > 0$ and $0 < c < 1$. The polynomials $\mathcal{M}_n(x;\beta, c)$, for $n= 0,1,2,\ldots$, satisfy the orthogonality relations:
\begin{equation}
 \sum_{k=0}^\infty \mathcal{M}_n(k;\beta,c) \mathcal{M}_m(k;\beta,c) (1-c)^\beta (\beta)_k \frac{c^k}{k!} = 0, \quad \text{for } n \neq m.\label{eq27}
\end{equation}
It is evident from (\ref{eq25}) and (\ref{eq26}) that when $d=1$ and $Q(z)=\frac{c}{c-1}z$, the polynomials $\mathcal{M}_m(k;\beta,c)$ and $\widehat{P}_n(k)$ are related as follows:
\[
\widehat{P}_n(k) = \widehat{\mathcal{M}}_n(k;\beta,c) = \frac{c^n(\beta)_n}{(c-1)^n} \mathcal{M}_n(k;\beta,c),
\]
where $\widehat{\mathcal{M}}_n(k)$ denotes the monic polynomials associated with $\mathcal{M}_n(k;\beta,c)$.
\subsection{Lowering operator}
By considering the matrix elements of the operator $SJ_-$, we have
\begin{equation}
\langle \beta,k|SJ_-|n,\beta\rangle=\sqrt{n(n+\beta-1)}\psi_{n-1,k}.\label{eq28}
\end{equation}
It is easy to see from \ref{relation2} that
\begin{align*}
SJ_-S^{-1}&=e^{ J_+}J_-e^{- J_+}\\
&=J_--2 J_0+J_+.
\end{align*}
Consequently,
\begin{align*}
\langle \beta,k|SJ_-|n,\beta\rangle&=\langle \beta,kn|(SJ_-S^{-1})S|n,\beta\rangle\\&=\langle \beta,k|J_+S|n,\beta\rangle+\langle \beta,k|J_-S|n,\beta\rangle-2\langle \beta,k|J_0S|n,\beta\rangle.
\end{align*}
Then we get
\begin{align*}
\langle \beta,k|SJ_-|n,\beta\rangle=\sqrt{(k+1)(k+\beta)}\psi_{n,k+1}-2(k+\frac{\beta}{2})\psi_{nk}+\sqrt{k(k+\beta-1)}\psi_{n,k-1}.\label{eq29}
\end{align*}
Put
\[
{\widehat{P}_n(x)} = \sqrt{\frac{n!(\beta)_n (\beta)_k}{k!}} \psi_{nk}.
\]
Then from \eqref{eq28} we have
\begin{align*}
n(n+\beta-1){\widehat{P}_{n-1}(k)} &= (k+\beta){\widehat{P}_n(k+1)} - (2k+\beta){\widehat{P}_n(k)} + k{\widehat{P}_n(k-1)} \\
&= kT_-\Delta^2 {\widehat{P}_n(k)} + \beta\Delta {\widehat{P}_n(k)},
\end{align*}
where $$\Delta f(x)=f(x+1)-f(x),\quad T_-f(x)=f(x-1).$$
So, the lowering operator for the polynomials \(\widehat{P}_n(k)\) is defined by the expression $$\sigma = k \Delta^2 T_- + \beta \Delta.$$ Consequently, it follows that:
\[
\sigma \widehat{P}_n(x) = n(n+\beta-1) \widehat{P}_{n-1}(x).
\]

\section{$d$-orthogonal polynomials of Meixner type}
In this section, we will interested with the determination of the explicit expression of a linear functional vector $({\mathcal{L}}_0,{\mathcal{L}_1},\,\dots,\,{\mathcal{L}}_{d-1})$ ensuring the $d$-orthogonality of the polynomial set $\mathcal{M}_n(x;\beta,c,d)$ generated by
 \begin{equation}
\sum_{n=0}^{\infty} \mathcal{M}_{n}(x;\beta,c,d)\frac{z^n}{n!}=e^{z+\frac{(-1)^r}{1-c}z^r}\hyp11{-x}{\beta}{-z}.\label{eq31}
\end{equation}
 $\mathcal{M}_n(x;\beta,c,d)$, which is reduced when $d=1$ to the Meixner polynomials of the first kind  is called $d$-OPS of Meixner type.\\
Now, we state our main result.
{
\begin{theorem}\label{thm1}The polynomial set $\mathcal{M}_n(x;\beta,c,d)$ generated by (\ref{eq31}) is $d$-OPS with respect to the linear functional vector given by
\begin{equation}
{\mathcal{L}}_i(f)=\sum_{k=0}^\infty w_{i\,k}\frac{(\beta)_k}{k!}f(k),\label{eq32}
\end{equation}
where
\begin{align*}
w_{i\,k}&=\frac{(-1)^i(1-c)^\frac{\beta+i}{r}}{r\Gamma(\beta)\sqrt{i!(\beta)_i}}
\frac{}{}\sum_{s=0}^{r-1}\frac{(-k)_s}{s!(\beta)_s}\Gamma\bigg(\frac{\beta+i+s}{r}\bigg)
\\\times &{}_{r+2}F_{2r}\left(
          \begin{array}{cc}
            1,\Delta(r,-k+s),\frac{\beta+s+i}{r}& \\
           \Delta(r,\beta+s),\Delta(r,s+1) &  \\
          \end{array}; \frac{(1-c)^{1/r}}{r^r}\right)
\end{align*}
and $$\Delta(r,m)=\Big(\frac{m}{r},\frac{m+1}{r},\ldots,\frac{m+r-1}{r}\Big).$$
\end{theorem}\label{thm1}
In particular, when  $r=d=1$ (then $s=i=0$), we have with the help of the binomial theorem that
\begin{align*}
w_{0k}&=(1-c)^\beta\hyp32{1,-k,\beta}{\beta,1}{1-c}\\
&=(1-c)^\beta\hyp10{-k}{-}{1-c}\\
&=(1-c)^\beta\sum_{n=0}^\infty\frac{(-k)_n}{n!}(1-c)^n\\
&=(1-c)^\beta c^k.
\end{align*}

Theorem (\ref{thm1}) is reduced to
$${\mathcal{L}_0}(f)=(1-c)^\beta\sum_{k=0}^\infty\frac{c^k(\beta)_k}{k!}f(k).$$
Hence we recognize the orthogonality of the Meixner polynomials.\\
In order to prove Theorem (\ref{thm1}), we need some preliminaries developments.
\subsection{Recurrence relation of $\phi_{nk}$ and linear functional vector}
We consider the matrix elements of the operator $S^{-1}J_0$. After similar calculus used in Section 4, we have
\begin{align*}
\langle \beta,n|S^{-1}J_0|k,\beta\rangle=(k+\frac{\beta}{2}){\phi}_{nk}.
\end{align*}
On the other hand,
\begin{align*}
\langle \beta,n|S^{-1}J_0|k,\beta\rangle&=\langle \beta,n|(S^{-1}J_0S)S^{-1}|n,\beta\rangle
\\&=\langle \beta,n|J_+S^{-1}|n,\beta\rangle+\langle \beta,n|J_0\sum_{i=0}^da_iJ_-^iS^{-1}|n,\beta\rangle\\&+\langle \beta,n|(\sum_{i=0}^da_iJ_-^i)S^{-1}|n,\beta\rangle.
\end{align*}
Therefore we arrive to the recurrence relation
\begin{equation}
(k+\frac{\beta}{2}){\phi}_{nk}=\sqrt{n(n+\beta-1)}{\phi}_{n-1,k}+\sum_{i=0}^{d}\alpha_{ni} {\phi}_{n+i,k},\label{eq33}
\end{equation}
where $\alpha_{ni}$ real numbers,\ $\alpha_{nd}\not=0$.\\
According to (\ref{eq33}), by the same manner as in \cite{V2}, we state the following:
\begin{proposition}  The coefficients $\phi_{nk}$
can be expressed as
\begin{equation}
\phi_{nk}=\sum_{i=0}^{d-1}\phi_{ik}R_n^{(i)}(k),\label{eq34}
\end{equation}
 where $R_n^{(i)}(k)$ are polynomials of argument $k$.\\
The degrees of these polynomials depend on $n$ in the following manner:\\
Assume that $n=dj+q$ where $q=0,\ldots,d-1$. Then
\begin{align*}
&\deg R_n^{(i)}=j\ \ \ \quad{if}\ \ i\leq q,\\
&\deg R_n^{(i)}=j-1\ \ \ \quad{if}\ \ i > q.
\end{align*}
\end{proposition}
\indent
The proof of this proposition proceeds by induction on $n$. Note that, it is sufficient to apply (\ref{eq33}) for $n=d$ to obtain
\begin{equation*}
\phi_{d,k}=\frac{1}{\alpha_{0d}}\big((k+\frac{\beta}{2})\phi_{0k}-\sum_{k=0}^{d-1}\alpha_{0i}\phi_{ik}\big)=\sum_{k=0}^{d-1}\phi_{ik}R_d^{(i)}(k),
\end{equation*}
where
\begin{equation*}
R_d^{(0)}(k)=\frac{k+\frac{\beta}{2}-\alpha_{0,0}}{\alpha_{0d}},\quad  R_d^{(i)}=-\alpha_{0i},\ i=1,...,d-1.
\end{equation*}
It is clear that $j=1$, $q=0$, $\deg R_d^{(0)}=1$,\quad $\deg R_d^{(i)}=0$.\\
Now, we consider the linear functional vector
 $({\mathcal{L}}_0,{\mathcal{L}_1}_,...,{\mathcal{L}}_{d-1})$  defined by
 \begin{equation}
 {\mathcal{L}_i}(f)=\sum_{k=0}^\infty\psi_{0k}\phi_{ik}f(k), \quad \text{with} \ i=1,...,d-1.\label{eq35}
\end{equation}
Then we have the following:
\begin{proposition}The $d$-OPS $\{\hat{P_n}(x)\}$ satisfies the following vector orthogonality relations:
\begin{align}
&{\mathcal{L}_i}(x^m\hat{P_n}(x))\neq0,\quad n=md+i,\label{eq36}\\
&{\mathcal{L}_i}(x^m\hat{P_n}(x))=0,\quad  n\geq md+i+1.\label{eq37}
\end{align}
\end{proposition}
In fact, relations (\ref{eq36}) and (\ref{eq37}) are immediate from (\ref{eq35}) and (\ref{eq34}).
\subsection{\bf{Explicit expression of  ${\bf{\phi}}_{ik}$}}
Now, we express explicitly the coefficients  $\phi_{ik}$ in terms of hypergeometric functions.\\
According to (\ref{eq23}), we have
\begin{align*}
\sum_{k=0}^\infty\sum_{n=0}^\infty\frac{\phi_{ik}\psi_{nk}}{\sqrt{n!(\beta)_n}}z^n&=e^{Q(z)}\sum_{k=0}^\infty
\sqrt{\frac{(\beta)_k}{k!}}\phi_{ik}\hyp11{-k}{\beta}{-z}\\
&=e^{Q(z)}\sum_{k=0}^\infty
\sqrt{\frac{k!}{(\beta)_k}}\phi_{ik}L_k^{\beta-1}(-z)
\end{align*}
where $L_k^{\beta-1}(z)$ are the Laguerre polynomials given by
 $$L_k^{\beta-1}(z)=\frac{(\beta)_k}{k!}\hyp11{-k}{\beta}{z}.$$
On the other hand, from the biorthogonality relation (\ref{eq13}), we have
\begin{align*}
\sum_{k=0}^\infty\sum_{n=0}^\infty\frac{\phi_{ik}\psi_{nk}}{\sqrt{n!(\beta)_n}}z^n=&
\sum_{n=0}^\infty\frac{z^n}{\sqrt{n!(\beta)_n}}\sum_{k=0}^\infty\phi_{ik}\psi_{nk}\\
=&\frac{z^i}{\sqrt{i!(\beta)_i}}.
\end{align*}
Then, we get
$$\sum_{k=0}^\infty\sqrt{\frac{k!}{(\beta)_k}}\phi_{ik}L_k^{\beta-1}(-z)=\frac{e^{-Q(z)}z^i}{\sqrt{i!(\beta)_i}}.$$
Also from the orthogonality of  $L_k^{\beta-1}(-z)$, the coefficients $\phi_{ik}$ are expressed as
$$\sqrt{\frac{k!}{(\beta)_k }}\phi_{ik}=\frac{(-1)^ik!}{\sqrt{i!(\beta)_i}\Gamma(\beta+k)}\int_0^\infty e^{-t-Q(-t)}t^{\beta+i-1}L_k^{\beta-1}(t)dt.$$
Using the facts that
\begin{align*}
&Q(t)=t+\frac{(-1)^r}{1-c}t^r\ ,\ L_k^{\beta-1}(t)=\frac{(\beta)_k}{k!}\sum_{m=0}^\infty\frac{(-k)_m}{m!(\beta)_m}t^m,\\
&\Gamma(\beta+k)=(\beta)_k\Gamma(\beta)\ ,\ x=\frac{t^r}{1-c}\\
& (mr+s)!=s!r^{mr}\prod_{j=0}^{r-1}\bigg(\frac{s+j}{r}\bigg)_m,\ \\
&(-k)_{mr+s}=(-k)_s!r^{mr}\prod_{j=0}^{r-1}\bigg(\frac{-k+s+j}{r}\bigg)_m,
\end{align*}
we directly obtain
\begin{align*}
\sqrt{\frac{k!}{(\beta)_k }}\phi_{ik}&=\frac{(-1)^ik!(\beta)_k}{\sqrt{i!(\beta)_i}k!\Gamma(\beta+k)}\sum_{m=0}^\infty\frac{(-k)_m}{m!(\beta)_m}
\int_0^\infty e^{\frac{-t^r}{1-c}}t^{\beta+m+i-1}dt\\&
=\frac{(-1)^i}{r\sqrt{i!(\beta)_i}\Gamma(\beta)}\sum_{m=0}^\infty\frac{(-k)_m}{m!(\beta)_m}(1-c)^\frac{\beta+m+i}{r}\Gamma\bigg(\frac{m+\beta+i}{r}\bigg)\\&
=\frac{(-1)^i(1-c)^\frac{\beta+i}{r}}{r\Gamma(\beta)\sqrt{i!(\beta)_i}}
\sum_{s=0}^{d-1}\sum_{m=0}^\infty\frac{(-k)_{mr+s}}{(\beta)_{mr+s}{(mr+s)}!}(1-c)^{m+\frac{s}{r}}\Gamma\bigg(m+\frac{\beta+i+s}{r}\bigg)\\&
=\frac{(-1)^i(1-c)^\frac{\beta+i}{r}}{r\Gamma(\beta)\sqrt{i!(\beta)_i}}
\sum_{s=0}^{r-1}\sum_{m=0}^\infty\frac{(-k)_s \big(\frac{\beta+i+s}{r}\big)_m\Gamma
\big(\frac{\beta+i+s}{r}\big)\prod_{j=0}^{r-1}\big(\frac{-k+s+j}{r}\big)(1-c)^{\frac{m}{r}}}
{s!(\beta)_sr^{mr}\prod_{j=0}^{r-1}\big(\frac{\beta+s+j}{r}\big)_m\prod_{j=0}^{r-1}\big(\frac{s+j+1}{r}\big)_m },
\end{align*}
Therefore
\begin{align*}
\psi_{0k}\phi_{ik}&=\sqrt{\frac{(\beta)_k}{k! }}\phi_{ik}\\
&=\frac{(-1)^i(1-c)^\frac{\beta+i}{r}(\beta)_k}{r\Gamma(\beta)\sqrt{i!(\beta)_i}k!}
\sum_{s=0}^{r-1}\frac{(-k)_s}{s!(\beta)_s}\Gamma\bigg(\frac{\beta+i+s}{r}\bigg)
\\\times &{}_{r+2}F_{2r}\left(
          \begin{array}{cc}
            1,\Delta(r,-k+s),\frac{\beta+s+i}{r}& \\
           \Delta(r,\beta+s),\Delta(r,s+1) &  \\
          \end{array}; \frac{(1-c)^{1/r}}{r^r}\right).
\end{align*}

\section*{Acknowledgment}
The authors would like to thank Professor Ali Zaghouani for his important remarks and suggestions that significantly improved this work.

\section*{Funding}This work is supported by  the ``Research Supporting Project number (RSPD2024R974), King Saud University, Riyadh, Saudi Arabia''.

\end{document}